\documentclass[12pt]{article}
\usepackage{amsfonts,amssymb,amsmath}
\usepackage{slashed}
\newcommand{\C}{\mathbb C}

\renewcommand{\H}{\mathbb H}
\newcommand{\N}{\mathbb N}
\renewcommand{\P}{\mathbb P}
\newcommand{\R}{\mathbb R}
\newcommand{\Z}{\mathbb Z}
\newcommand{\Aut}{\mbox{Aut}}
\newcommand{\End}{\mbox{End}}
\newcommand\mb[1]{\mbox{\rm#1}}
\newcommand{\tr}{\mbox{tr}}
\newcommand{\id}{\mathbb I}

\renewcommand{\)}{\mb{)}}
\newtheorem{definition}{Definition}[subsection]
\newtheorem{lemma}[definition]{Lemma}
\newtheorem{prop}[definition]{Proposition}
\newtheorem{theorem}[definition]{Theorem}
\newenvironment{proofsketch}{\textsl{Proof \(sketch\):\hspace*{\fill}\\ {}}}{\hspace*{\fill}$\Box$\\[0.5em]}
\numberwithin{equation}{subsection}
\begin{document}
\title{K3 en route\\ From Geometry to Conformal Field Theory}
\author{Katrin Wendland}
\date{}
\maketitle
\begin{abstract}
To pave the way for the journey from geometry to
conformal field theory (CFT), these notes present the background 
for some basic CFT constructions from
Calabi-Yau geometry. Topics include the complex 
and K\"ahler geometry of Calabi-Yau manifolds and 
their classification in low dimensions. 
I furthermore discuss CFT 
constructions for the simplest known examples 
that are based in Calabi-Yau geometry, namely 
for the toroidal superconformal field theories and 
their $\Z_2$-orbifolds.
En route from geometry to CFT, I 
offer a discussion of K3 surfaces as the simplest class
of Calabi-Yau manifolds where non-linear sigma model
constructions bear mysteries to the very day.
The elliptic genus in CFT and in 
geometry is recalled as an instructional
piece of evidence in favor of a deep connection
between geometry and conformal field theory.
\end{abstract}
\section*{Introduction}
These lecture notes aim to make a contribution to paving the way
from geometry to conformal field theory.\\[0.5em]
While two-dimensional conformal field theories can 
be defined abstractly in mathematics without 
reference to the algebraic geometry of complex
manifolds of dimension greater than one, some of
the most interesting examples are expected to arise
geometrically from Calabi-Yau manifolds
of complex dimension two or higher.
Indeed, Calabi-Yau manifolds
are the candidates for consistent geometric backgrounds for superstrings.
Roughly, the string dynamics in such a geometric background
are governed by so-called 
\textsc{non-linear sigma models}, whose equations of motion 
imply conformal invariance. The resulting
quantum field theory on the world-sheet is thus a 
superconformal field theory. 
From a mathematical point of view, this looks encouraging, since 
conformal field theory
allows an independent mathematical approach, while
string theory as a whole is not comprehensible to a 
mathematically satisfactory degree, as yet. 
However, non-linear sigma model constructions are
far from well understood, mathematically.
The only examples of (compact) Calabi-Yau manifolds
where explicit constructions of a non-linear
sigma model are known are the
complex tori and their orbifolds. 
It seems that we are jumping out 
of the frying pan into the fire:
By focusing on conformal field theory instead
of string theory, we trade one mathematically dissatisfactory
approach by a concept whose abstract definition prevents it from
being immediately applied. 
Though string theory in its full generality seems to allow
non-geometric phases, it is still crucially connected to geometry. 
To advance the subject, I am convinced that
it is a \textsl{Sine qua non} to get a better
understanding of the precise relation between 
geometry and conformal field theory, beginning with those cases
where constructions are explicitly known.\\[0.5em]
These lecture notes should be viewed as an invitation
to this journey from geometry to conformal field theory.
They give a lightning introduction to the subject,
as do many other excellent sources, so I try to
keep the exposition somewhat 
complementary to existing works. 
As key examples, which are
both sufficiently simple and mysterious,
K3 surfaces play a special
role en route from geometry to conformal field
theory.\\

\noindent
These notes are structured as follows:\\[0.5em]
Section \ref{CY} is devoted to some background
in Calabi-Yau geometry. 
In Section \ref{CYcdt}, 
I recall basic mathematical concepts leading
to the definition of Calabi-Yau manifolds, and some 
of their fundamental properties.
Section \ref{classification} introduces a number of 
topological  invariants of Calabi-Yau manifolds and
culminates in a summary of  mathematical arguments 
that yield the classification of
Calabi-Yau manifolds of complex dimension $D\leq2$.
The classification naturally leads to the
definition of K3 surfaces. 
An important 
topic whose discussion would naturally follow is the structure
of the  moduli spaces of complex structures, K\"ahler structures
and Hyperk\"ahler structures on K3 surfaces. However, 
since these topics have already been discussed elsewhere,
even by myself \cite{nawe00,we03},
I omit them in these notes.
Instead, Section \ref{Kummer} 
gives a detailed summary of the Kummer construction
as a classical example of a geometric orbifold procedure.
In particular,  I argue that the Kummer surfaces
constitute a class of K3 surfaces which are very well
understood,
because their geometric properties are entirely governed by 
much simpler
manifolds, namely the underlying  complex two-tori.

Section \ref{CFT} discusses explicit examples
of superconformal field theories that are obtained
from non-linear sigma model constructions. 
Since I have done so elsewhere, in these notes I do not
offer a proposal for the defining mathematical properties of 
conformal field theory. Instead,
the recent review \cite{we14} should be viewed
as the conformal field theory companion to these lecture notes,
while the present exposition, in turn, can be understood as the geometric companion
of \cite{we14}. 
Concretely, Section \ref{toroidal} gives a brief review of toroidal 
superconformal field theories. Lifting the Kummer construction
of Section \ref{Kummer} 
to the level of conformal field theory, Section \ref{Z2} addresses 
$\Z_2$-orbifold constructions of toroidal superconformal field 
theories. To substantiate the expectation that these orbifold
conformal field theories are correctly interpreted as non-linear
sigma models on Kummer surfaces, I 
include a brief discussion of elliptic genera:
I introduce the conformal field theoretic 
elliptic genus as a counter part of the geometric elliptic genus.
I then show that the geometric elliptic genus of K3 surfaces
agrees with the conformal field theoretic
elliptic genus of the $\Z_2$-orbifold 
conformal field theory obtained
from a toroidal superconformal field
theory on a complex two-torus, recalling the known
proof \cite{eoty89}.

The final Section \ref{K3theories} places K3
en route from geometry to conformal field theory:
I motivate and discuss the definition of K3 theories, which is 
formulated purely within representation theory. While this may be
mathematically satisfying, it entails a consideration of K3 theories 
against non-linear sigma models on K3, for which explicit direct
constructions on smooth K3 surfaces are lacking.
I recall the role of the chiral de Rham complex in the context 
of elliptic genera of Calabi-Yau manifolds and conformal field theories,
respectively. 
I conclude with a few speculations on the vertex algebra which can 
be obtained from the cohomology of the chiral de Rham complex,
as candidate for a recovery of some of the conformal field theory structure
from purely geometric ingredients.
\\[1em]
{\textbf{Acknowledgements}}\nopagebreak

My sincere thanks go to the
organizers Alexander Cardona, Sylvie Paycha, Andr\'es Reyes, Hern\'an Ocampo
and the participants
of the tremendously inspiring
8th Summer School on ``Geometric, Algebraic and Topological Methods
for Quantum Field Theory 2013'' at
Villa de Leyva, Columbia, where the lectures leading to these
notes were presented. 

I am grateful to the organizers of the 2nd  COST MP2010
Meeting ``The String Theory Universe'', 
particularly to the local organizer Gabriele Honecker, 
for inviting me to this stimulating event at the Mainz Institute for 
Theoretical Physics
(MITP). 
My thanks go to the MITP
for its hospitality and its partial support during the completion
of these lecture notes. 

This work was in part  supported by my 
ERC Starting Independent Researcher Grant StG No. 204757-TQFT
``The geometry of topological quantum field theories''.
\section{Calabi-Yau geometry}\label{CY}
By definition, a \textsc{Calabi-Yau manifold} is a compact\footnote{By \textsc{compact}, I mean bounded 
and closed; one also frequently finds definitions of non-compact Calabi-Yau manifolds, but those
will not be of interest in these lectures.}  K\"ahler manifold with trivial canonical bundle. 
This section gives an introduction to the mathematical ingredients of this definition and illustrates its consequences. 

Examples of Calabi-Yau manifolds are the complex tori, which also furnish the only examples where a non-linear sigma model construction immediately yields an associated conformal field theory (see Section \ref{toroidal}). In Section \ref{classification}, we will see that in complex dimension one and two, apart from tori, the only  Calabi-Yau manifolds are the \textsc{K3 surfaces}. The final subsection is devoted to the \textsc{Kummer construction}, yielding a special class of K3 surfaces which are almost as well under control from a conformal field theoretic point of view as are the complex tori (see Section \ref{Z2}).

There are a number of excellent books on the topics presented in this section, see for example
\cite{bhpv84, grha78,hu05,jo00,we73}.
In particular, if not stated otherwise, then  proofs of the classical results below
can be found in these references.
\subsection{The Calabi-Yau condition}\label{CYcdt}
In the following, familiarity with the concepts of differential geometry and 
complex analysis is assumed. The definition of topological and Riemannian
manifolds over $\R$ and of vector bundles is crucial and can be found in textbooks
like \cite{do92,le09,mo01,on83}. The concept of holomorphic functions
and their special properties is the foundation of complex analysis, where
good textbooks include \cite{frbu05,la99,re91}.

Recall that a \textsc{$D$-dimensional complex manifold} is a differentiable real  
$2D$-dimensional manifold $Y$  together with a holomorphic atlas 
$(U_\alpha, \varphi_\alpha)_{\alpha\in A}$, i.e.\ an open covering 
$\{U_\alpha\mid\alpha\in A\}$ of $Y$ with diffeomorphisms 
$\varphi_\alpha\colon U_\alpha\longrightarrow V_\alpha$, $V_\alpha\subset\C^D$ 
open for all $\alpha\in A$, such that all coordinate changes 
$\varphi_\alpha\circ\varphi_\beta^{-1}\colon
\varphi_\beta(U_\alpha\cap U_\beta)\longrightarrow \varphi_\alpha(U_\alpha\cap U_\beta)$ 
are holomorphic. 
Such a holomorphic atlas is said to define a \textsc{complex structure} on the real
manifold $Y$. There is an alternative description, which
sometimes is more convenient, and which is motivated
by the following observations:

Consider $y\in U_\alpha$ and
note that due to $\varphi_\alpha(U_\alpha)=V_\alpha\subset\C^D$,
for the tangent space $T_{\varphi_\alpha(y)} V_\alpha$
we have a natural identification
$\C^D=T_{\varphi_\alpha(y)} V_\alpha$, such that
multiplication by $i=\sqrt{-1}$ on $\C^D$ yields an endomorphism of $T_{\varphi_\alpha(y)} V_\alpha$.
On the tangent space $T_yY$ this induces 
an endomorphism $I\in\End(T_yY)$ with $D\varphi_\alpha(Iv)=iD\varphi_\alpha(v)\in\C^D$
for all $v\in T_yY$. The endomorphism $I$
is checked to be independent of the choice of (holomorphic!) coordinates. 
Its $\C$-linear extension to 
the complexification $T^\C Y := TY\otimes_\R \C$ of the tangent bundle $TY$
is also denoted by $I$.
By construction, $I^2=-\id$, such that we have an eigenspace 
decomposition $T^\C Y =T^{1,0}Y\oplus T^{0,1}Y$ for $I$, where $I$ acts on 
the fibers of $T^{1,0}Y$ by multiplication by $i$, and on those of 
$T^{0,1}Y$ by multiplication by $-i$. In fact, $T^{1,0}Y$ is a holomorphic 
vector bundle whose cocycles can be given by the Jacobians of 
the coordinate changes in any 
holomorphic atlas of $Y$, and $T^{0,1}Y=\overline{T^{1,0}Y}$
with the complex conjugation $\overline{v\otimes\lambda}:=v\otimes\overline\lambda$
for $y\in Y,\, v\in T_yY,\, \lambda\in\C$.
The induced decomposition  of the cotangent bundle
yields a decomposition $d=\partial+\overline\partial$ of the exterior
differential $d$, where on 
$\mathcal A^{p,q}(Y)$,
the space of $(p,q)$-forms on $Y$, we have 
$$
\partial\colon \mathcal A^{p,q}(Y)\longrightarrow
\mathcal A^{p+1,q}(Y)
\quad\mbox{ and }\quad
\overline\partial\colon \mathcal A^{p,q}(Y)\longrightarrow
\mathcal A^{p,q+1}(Y).
$$
Since $T^{1,0}Y$ is a holomorphic vector bundle,  it is closed under the Lie bracket
(see e.g.\ \cite[\S1 and App.\ B]{on83} for a discussion of the Lie bracket).

Vice versa, on a differentiable real  $2D$-dimensional manifold $Y$, a
fiber-wise endomorphism
$I\in\End(TY)$ with $I^2=-\id$ is called an \textsc{almost complex structure}. 
Given an almost complex structure, the fiber-wise ei\-gen\-space decomposition 
$T^\C Y =T^{1,0}Y\oplus T^{0,1}Y$ for $I$ holds
as above, with complex vector bundles $T^\C Y$, $T^{1,0}Y$ and $T^{0,1}Y$. 
The almost complex structure $I$ is called
\textsc{integrable} if the bundle $T^{1,0}Y$ closes under the Lie bracket. This
condition is equivalent to $d=\partial+\overline\partial$ on each
$\mathcal A^{p,q}(Y)$ as above, where on $\mathcal A^{p,q}(Y)$,
the operators $\partial,\,\overline\partial$ are obtained by composing 
the operator $d$ with the projector to $\mathcal A^{p+1,q}(Y)$ and $\mathcal A^{p,q+1}(Y)$,
respectively.
By the celebrated Newlander-Nirenberg theorem \cite{neni75}, every
integrable almost complex structure on $Y$ is induced by a unique
complex structure.

In summary, the choice of a complex structure on $Y$ is tantamount to
the choice of an endomorphism $I$ which defines an integrable almost 
complex structure on $Y$. \\[0.5em]
Next, Riemannian geometry comes into play: We choose a 
Riemannian metric $g$ on our complex manifold $Y$, 
viewed as a real manifold, and discuss meaningful additional
compatibility conditions.
Let $g$ also denote  the sesquilinear continuation of $g$ to $T^\C Y$,
where I use the convention 
$$
\forall y\in Y,\;\; u,\, v\in T_y^\C Y,\;\;  \lambda,\,\mu\in\C\colon\quad
g(\lambda u, \mu v)=\lambda\overline\mu g(u,v).
$$
With respect to local complex coordinates $(z^1,\ldots,z^D)$ on $U\subset Y$,
and restricting to real tangent vectors,
\begin{eqnarray*}
g_{\mid TU\times TU}
&=& {\textstyle{1\over2}} \sum_{j,k} \left(
g_{j,\overline k} dz^j\otimes d\overline z^k
+ g_{k,\overline\jmath}  d\overline z^j\otimes dz^k
 \right.\\[-1em]
&&\quad\quad\quad\quad\left.
+ g_{j,k}  dz^j\otimes dz^k
+ \overline{g_{j,k}}  d\overline z^j\otimes d\overline z^k
\right)
,\\[0.8em]
&&\quad\quad\mbox{ where }
{\textstyle g_{j,\overline k}:= 2g\left( {\partial\over\partial z^j}, {\partial\over\partial z^k}\right)},
\;\;
{\textstyle g_{j, k}:= 2g\left( {\partial\over\partial z^j}, {\partial\over\partial\overline z^k}\right)},
\end{eqnarray*}
with a Hermitean matrix $(g_{j,\overline k})_{j,k\in\{1,\ldots, D\}}$. The complex structure on $Y$
is said to be \textsc{compatible} with the metric $g$, if  the 
corresponding almost complex structure
$I\in\End(TY)$ is orthogonal with respect to $g$, in other words if
$$
\forall  y\in Y,\;\; u,\, v\in T_y Y\colon\quad
g(I u, I v)= g(u,v).
$$ 
Then $g_{j,k}=0$ 
for all $j,k\in\{1,\ldots, D\}$, so
$$
g_{\mid TU\times TU}
= {\textstyle{1\over2}} \sum_{j,k} \left(
g_{j,\overline k} dz^j\otimes d\overline z^k
+ g_{k,\overline\jmath}  d\overline z^j\otimes dz^k\right).
$$
One checks directly  that this compatibility condition  implies that
with respect to $g$, we have $T^{1,0}_yY\perp T^{0,1}_yY$ for all $y\in Y$.

If $I$ is compatible with  $g$, then the $\C$-linear extension $\omega$ 
of the bilinear form 
$$
\forall  y\in Y,\;\; u,\, v\in T_y Y\colon\quad
\omega(u,v):=g(Iu,v)
$$
to $T^\C Y$ is called the \textsc{K\"ahler form}. By 
construction, $\omega$ is real, and by a direct calculation
one checks that compatibility of $g$ with $I$ implies
that $\omega$ is antisymmetric. Denoting by 
$\mathcal A^k_\R(Y)$
the space of real $k$-forms on $Y$,
and with respect to local complex coordinates $(z^1,\ldots,z^D)$ on $U\subset Y$ as above,
we have
$$
\omega\in\mathcal A^{1,1}(Y)\cap\mathcal A^2_\R(Y),\quad
\omega_{\mid U}
= {\textstyle{i\over2}} \sum_{j,k} g_{j,\overline k} dz^j\wedge d\overline z^k.
$$
It is important to keep in mind that $g$ is extended to $T^\C Y$ as a \textsc{Hermitean
sesquilinear form}, while $\omega$ is \textsc{antisymmetric and $\C$-bilinear} on $T^\C Y$.
Note that any two ingredients of the triple $(I,\, g,\, \omega)$ determine the third one.\\[0.5em]
The following additional condition on the metric turns out to have far-reaching consequences \cite{ka33}:
\begin{definition}\label{Kc}
Consider a Riemannian manifold $(Y,g)$, equipped with a compatible
complex structure. Then the  metric $g$ is called 
\textsc{K\"ahler metric} if and only if the K\"ahler form $\omega$ is a closed differential form:
$d\omega=0$. 

A complex manifold $Y$ is called \textsc{K\"ahler manifold} if
a K\"ahler metric exists on $Y$.
\end{definition}
A harbinger of the fact that this condition is a very natural and 
interesting one is the following lemma, which states as many as
six equivalent formulations of the K\"ahler condition, each
emphasizing a different geometric aspect:
\begin{lemma}\label{Kaehlerequivalences}
Consider a Riemannian manifold $(Y,g)$ with compatible
complex structure $I$, such that $Y$ is a complex $D$-dimensional
manifold and $\omega$ its K\"ahler form. 
Then the following are equivalent:
\begin{enumerate}
\item
The metric $g$ is a K\"ahler metric.
\item
For every $y\in Y$, there is an open neighborhood $U\subset Y$
of $y$ and a smooth function $f\colon U\rightarrow\R$, such that 
$\omega_{\mid U}=i\partial\overline\partial f$.
\item
The Levi-Civita connection for $g$, that is, the unique 
torsion-free metric connection, agrees with the 
Chern \(or holomorphic\) connection for $I$, i.e.\ with the unique metric
connection whose $(0,1)$ part is $\overline\partial$.
\item 
The almost complex structure $I$ of $Y$ is parallel
with respect to the Levi-Civita connection for $g$.
\item
With respect to arbitrary holomorphic coordinates
$(z^1,\ldots, z^D)$ on $U\subset Y$, the coefficients of $\omega_{\mid U}
= {i\over2} \sum_{j,k} g_{j,\overline k} dz^j\wedge d\overline z^k$ 
obey
${\partial g_{k,\overline l}\over\partial z^m} = {\partial g_{m,\overline l}\over\partial z^k}$
for all $k,\,l,\,m\in\{1,\ldots,D\}$,
or equivalently
${\partial g_{k,\overline l}\over\partial\overline z^m} 
= {\partial g_{k,\overline m}\over\partial\overline z^l}$
for all $k,\,l,\,m\in\{1,\ldots,D\}$.
\item
The metric \textsc{osculates to second order} with the
standard Euclidean metric, that is, for every $y\in Y$, there are
holomorphic coordinates $(z^1,\ldots, z^D)$
around $y$ such that $z(y)=(z^1,\ldots, z^D)(y)=0$ and
$$
\textstyle
g\left( {\partial\over\partial z^j}, {\partial\over\partial  z^k}\right)
=\delta_{j,k} + {\mathcal O}(|z|^2) \quad \forall j,\, k\in\{1,\ldots,D\}.
$$
\item
The holonomy representation of $Y$ is unitary on each tangent space.
\end{enumerate}
\end{lemma}
In short, K\"ahler metrics are very similar to the Euclidean metric.\\[0.5em]
The K\"ahler condition has a number of important consequences.
For example,
the standard Laplace operators on a K\"ahler manifold 
obey a very simple relation:
$\Delta_d=2\Delta_\partial=2\Delta_{\overline\partial}$, where
$\Delta_d=dd^\ast+d^\ast d$, 
$\Delta_{\partial}=\partial\partial^\ast+\partial^\ast\partial$, and
$\Delta_{\overline\partial}
=\overline\partial\,\overline\partial^\ast+\overline\partial^\ast\overline\partial$.
Here, $d^\ast$, $\partial^\ast$ and $\overline\partial^\ast$ can be viewed
as the $L^2$-duals of $d$, $\partial$ and $\overline\partial$ if $Y$ is
compact.
One furthermore has
\begin{lemma}\label{Kaehlerconsequences}
Consider a compact complex K\"ahler manifold  $Y$, and denote its
de Rham and Dolbeault cohomology groups, respectively, by
\begin{eqnarray*}
H^k(Y,\C) 
&\!\!\!=\!\!\!& \mbox{\rm ker}(d:\mathcal A^k(Y)\longrightarrow\mathcal A^{k+1}(Y))\slash
\mbox{\rm im}(d:\mathcal A^{k-1}(Y)\longrightarrow\mathcal A^{k}(Y)),\\
H^{p,q}(Y,\C) 
&\!\!\!=\!\!\!& \mbox{\rm ker}(\overline\partial:\mathcal A^{p,q}(Y)\longrightarrow\mathcal A^{p,q+1}(Y))\slash
\mbox{\rm im}(\overline\partial:\mathcal A^{p,q-1}(Y)\longrightarrow\mathcal A^{p,q}(Y))
\end{eqnarray*}
with $\mathcal A^{k}(Y)$
the space of $k$-forms on $Y$. 
Then for all $k,\,p,\,q\in\N$,
$$
H^k(Y,\C)  = \bigoplus_{r+s=k} H^{r,s}(Y,\C), \quad
H^{p,q}(Y,\C) = \overline{H^{q,p}(Y,\C)}.
$$
\end{lemma}
Apart from the Euclidean spaces and complex tori, important examples of K\"ahler
manifolds are the complex projective spaces $\P^D$. Indeed, the
\textsc{Fubini-Study metric} defines a K\"ahler metric on $\P^D$; with respect
to homogeneous coordinates $(z^0\colon\cdots\colon z^D)$, its
K\"ahler form on the standard coordinate neighborhood
$U_j:=\{ (z^0\colon\cdots\colon z^D)\in\P^D \mid z^j\neq0 \}$ 
is
\begin{equation}\label{FS}
\textstyle
\omega_{\mid U_j}={i\over2\pi} \partial\overline\partial 
\log\left( \smash{\sum\limits_{k=0}^D \left| {z^k\over z^j}\right|^2} 
\vphantom{\left| {z^k\over z^j}\right|} \right).
\end{equation}
Since the restriction of a K\"ahler metric on a complex manifold
$Y$ to a complex submanifold yields a K\"ahler metric on the submanifold, this implies that
all closed algebraic manifolds are compact K\"ahler manifolds.\\[0.5em]
\noindent
Finally, we are ready to state
\begin{definition}\label{CYdef}
A  manifold $Y$ is called \textsc{Calabi-Yau manifold} if 
$Y$ is a compact K\"ahler manifold with trivial canonical bundle $K_Y$.
\end{definition}
Recall that the canonical bundle $K_Y$ of a complex $D$-manifold $Y$ is
the determinant  line bundle of the holomorphic cotangent bundle $\Omega_Y:=(T^{1,0}Y)^\ast$,
i.e.\ $K_Y=\Lambda^D\Omega_Y$. Hence a  K\"ahler manifold is Calabi-Yau
if and only if $h^{D,0}(Y):=\dim_\C H^{D,0}(Y,\C)=1$, or equivalently if the holonomy
representation of $Y$ is special unitary on each tangent space. 

Examples of Calabi-Yau manifolds include all complex tori $\C^D/L$
($L\subset\C^D$ a lattice of rank $2D$) and the degree $D+2$
hypersurfaces in $\P^{D+1}$.\\[0.5em]
It is  natural to investigate whether a given Calabi-Yau manifold
$Y$ possesses any K\"ahler metrics with particularly appealing properties. 
Indeed, the search within \textsc{K\"ahler classes} turns out to be fruitful.
Here, I make use of the fact that the K\"ahler form $\omega$ of a 
K\"ahler metric on $Y$ is a real, closed $(1,1)$-form. Thus, $\omega$
represents a class  $[\omega]\in H^{1,1}(Y,\C)\cap H^2(Y,\R)$, which in fact is non-zero if $Y$
is compact. Another K\"ahler metric on $Y$ is said to \textsc{belong to 
the same K\"ahler class as $g$}, if its K\"ahler form represents the same
class $[\omega]\in H^{1,1}(Y,\C)\cap H^2(Y,\R)$.
One now has the seminal 
\begin{theorem}[Calabi-Yau Theorem \cite{ca54,ya78}]\label{CYthm}
\hspace*{\fill}Consider \hspace*{\fill}a\linebreak Ca\-la\-bi-Yau manifold $Y$ with K\"ahler metric $g$. 
Then there exists a unique \textsc{Ricci-flat K\"ahler metric}
 in the K\"ahler class of $g$.
 \end{theorem}
The proof of this theorem is non-constructive \cite{ya78}; 
with the exception of complex tori, explicit forms of Ricci-flat K\"ahler metrics
on smooth  Calabi-Yau manifolds are not known.
\subsection{Classifying Calabi-Yau manifolds}\label{classification}
A classification of Calabi-Yau $D$-manifolds turns out to be possible in dimension
$D\leq2$, as I shall explain in this section.
In the following, if not stated otherwise,
let $Y$ denote a connected Calabi-Yau $D$-manifold.
I first recall some basic topological invariants of $Y$.

The theory of elliptic differential operators ensures that the cohomology
of $Y$ is finite dimensional,
$h^{p,q}(Y):=\dim_\C H^{p,q}(Y,\C)<\infty$ for all $p,\,q\in\N$.
These so-called \textsc{Hodge numbers} are topological invariants
which enjoy a number of constraints for our connected Calabi-Yau $D$-manifold $Y$:
\begin{equation}\label{hodgenumbers}
\begin{array}{rccccl}
h^{0,0}(Y) &=& h^{D,D}(Y) &=& h^{D,0}(Y) &=\;h^{0,D}(Y)\;=\;1,\quad\\[0.5em]
h^{p,q}(Y) &=& h^{q,p}(Y) &=& h^{D-p,D-q}(Y) &\quad\forall p,\,q\in\{0,\ldots,D\}.
\end{array}
\end{equation}
Here, $h^{p,q}(Y) = h^{q,p}(Y)$ is an immediate consequence of Lemma \ref{Kaehlerconsequences},
while $h^{p,q}(Y) = h^{D-p,D-q}(Y)$ is a consequence of the \textsc{Serre duality} \cite{se55}. 
$h^{0,0}(Y)=h^{D,D}(Y)=1$ is true by assumption, since $Y$ is compact and connected, while 
$h^{D,0}(Y)=h^{0,D}(Y)=1$ follows from  the triviality of the canonical bundle.

Using the Hodge numbers, one obtains the following classical topological 
invariants:
\begin{definition}\label{Eulerdef}
For a compact K\"ahler manifold $Y$ of dimension $D$,  the
\textsc{Euler characteristic} $\chi(Y)$, 
\textsc{holomorphic Euler characteristic} $\chi(\mathcal O_Y)$
and \textsc{signature} $\sigma(Y)$ are defined by
$$
\begin{array}{rclrcl}
\chi(Y) &:=& \sum\limits_{p,q=0}^D (-1)^{p+q} h^{p,q}(Y),\quad
&\chi(\mathcal O_Y) &:=& \sum\limits_{q=0}^D (-1)^{q} h^{0,q}(Y),\quad\\[1em]
&&&\sigma(Y) &:=& \sum\limits_{p,q=0}^D (-1)^{q} h^{p,q}(Y).
\end{array}
$$
For any holomorphic vector bundle $E\rightarrow Y$, 
$$
\chi(E):=\sum_{q=0}^D (-1)^q \dim H^q(Y,E)
$$ 
is the \textsc{holomorphic Euler characteristic} of $E$.
\end{definition}
By (\ref{hodgenumbers}), the signature $\sigma(Y)$ vanishes if the complex dimension
$D$ of $Y$ is odd; I have thus trivially extended the traditional definition of the signature of
oriented compact manifolds whose real dimension is divisible by $4$ to all compact
complex K\"ahler manifolds. 

We are now in a position to state a first classification result:
\begin{theorem}\label{CYonefolds}
Every connected Calabi-Yau one-manifold is biholomorphic to a torus
$\C/L$ with $L\subset\C$ a lattice of rank $2$.
\end{theorem}
\begin{proofsketch}
Since $D=1$, the constraints (\ref{hodgenumbers}) already fix all Hodge numbers.
Hence by Definition \ref{Eulerdef} we have $\chi(Y)=0$. 
Therefore, the claim follows from the  classification of compact Riemann surfaces.
A discussion of this deep and fundamental
classification result can be found in textbooks on differential
topology, for example \cite{fo81,hi76}. In particular, one needs to use the fact that
$\chi(Y)$ agrees with the Euler-Poincar\'e characteristic of $Y$ known from topology.
This in turn 
is deeply linked to the relation between de Rham cohomology, \u Cech cohomology
and singular homology, which is explained in depth in \cite{mo01}.
\end{proofsketch}
Another type of topological invariants 
are the \textsc{characteristic classes}, all of which 
can be traced back to \textsc{Chern classes} 
for our  Calabi-Yau manifold $Y$.
Recall that  by definition, a complex vector bundle 
$E\rightarrow Y$ of rank $r$
has Chern classes $c_k(E)\in H^{2k}(Y,\Z)$, $k\in\{0,\ldots,D\}$,
which are uniquely determined by the following four conditions (a)-(d):
\begin{enumerate}\renewcommand{\theenumi}{(\alph{enumi})}
\item
$c_0(E)=[1]$.  
\item
For the dual $\mathcal O(1)\rightarrow \P^D$ of the 
tautological line bundle, $c_1(\mathcal O(1))$ is the K\"ahler class 
of the Fubini-Study metric on $\P^D$ with K\"ahler form (\ref{FS}).
\item 
For
smooth $f\colon X\longrightarrow Y$ and  $k\in\{0,\ldots,D\}$, one has
$c_k(f^\ast E)=f^\ast c_k(E)$. 
\item
The
\textsc{total Chern class}
\smash{$c(E):=\sum\limits_{k=0}^D c_k(E)=\sum\limits_{k=0}^r c_k(E)$}
for line bundles $L_1,\ldots, L_r$ on $Y$ obeys
$$
c(L_1\oplus\cdots\oplus L_r) = c(L_1)\wedge\cdots\wedge c(L_r).
$$
\end{enumerate}
It follows that the total Chern class of every trivial bundle is $[1]\in H^0(Y,\Z)$.
By convention, the Chern classes $c_k(Y)\in H^{2k}(Y,\Z)$ of a complex manifold $Y$ are the
Chern classes of its holomorphic tangent bundle $T:=T^{1,0}Y$. 
One checks that $c_1(Y)=-c_1(K_Y)$, where as above,
$K_Y$ is the canonical bundle of $Y$.
Hence every Calabi-Yau manifold  $Y$ has vanishing first Chern class $c_1(Y)=0$. 
Vice versa, using the exponential sheaf sequence one finds:
If $Y$ is a compact K\"ahler manifold with $c_1(Y)=0$ and vanishing first
Betti number, $b_1(Y)=\dim_\C H^1(Y,\C)=0$, then $Y$ is Calabi-Yau\footnote{Instead
of Definition \ref{CYdef}, some authors
give a more restrictive definition of Calabi-Yau manifolds: 
A complex $D$-manifold $Y$ is then called Calabi-Yau, if and only if $Y$ is a 
compact K\"ahler manifold
with $h^{1,0}(Y)=0$ and $c_1(Y)=0$ or, equivalently, such that the holonomy representation
on every tangent space is an irreducible representation of $SU(D)$.
Hence complex tori are not Calabi-Yau according to this definition,
whereas I view them as the simplest examples of Calabi-Yau manifolds.}.

If $E\rightarrow Y$ is a complex vector bundle of rank $r$, then
by the so-called \textsc{splitting principle}, we may work with \textsc{formal Chern
roots} $e_1,\ldots,e_r$, such that $c(E)=\prod_{j=1}^r(1+e_j)$, where the $e_j$ 
are elements of a ring extension of $H^\ast(Y,\R)$. 
It then follows that $c_k(Y)=\sigma_k(e_1,\ldots,e_r)$
with the elementary symmetric polynomials $\sigma_k$. 
More generally, for an analytic function $f$ on 
a neighborhood of $0$ in $\C^r$, by $f(e_1,\ldots,e_r)$ one
denotes the power series expansion of $f$ about the origin with insertions
$(e_1,\ldots,e_r)$.
This allows the definition of further  topological invariants:
\begin{definition}\label{cherninv}
Consider a compact complex $D$-manifold $Y$ with holomorphic tangent
bundle $T:=T^{1,0}Y$, and a complex vector bundle
$E\rightarrow Y$ of rank $r$. Let $y_1,\ldots,y_D$ and $e_1,\ldots,e_r$ denote
the formal Chern roots, such that $c(Y)=\prod_{j=1}^D(1+y_j)$
and $c(E)=\prod_{j=1}^r(1+e_j)$. 

\noindent
Then the \textsc{Todd genus} of $Y$ is given by
$$
\mbox{Td}(Y) 
:= \prod_{j=1}^D {y_j\over 1-\exp(-y_j)}
= c_0(Y)+{\textstyle{1\over2}} c_1(Y) + {\textstyle{1\over12}} (c_1(Y)^2+c_2(Y)) + \cdots
$$ 
The \textsc{Chern character} of the bundle $E$ is
$$
\mbox{ch}(E) := \sum_{j=1}^r \exp(e_j) 
= r\,c_0(Y)+ c_1(E) + {\textstyle{1\over2}} (c_1(E)^2-2c_2(E)) + \cdots
$$
\end{definition}
On first sight, the topological invariants that are obtained from the Hodge numbers 
in Definition \ref{Eulerdef} and those that are obtained from the Chern classes
in Definition \ref{cherninv} have very different flavors: While the former are integers,
the latter are  cohomology classes. However, evaluation of such
cohomology classes
on the fundamental cycle
of $Y$, denoted by $\int_Y$, yields integers from integral cohomology classes.
Here, it is understood 
that $\int_Y \alpha=0$ if $\alpha\in H^k(Y,\Z)$ with $k\neq2D$.
The \textsc{Atiyah-Singer Index Theorem} \cite{atsi63} governs the deep relationship between
the two types of invariants; in the context of complex manifolds and holomorphic
bundles, which is relevant to our discussion, 
a precursor of this theorem is the following seminal
\begin{theorem}[Hirzebruch-Riemann-Roch Formula \cite{hi54}]\label{ashrr}\hspace*{\fill}\\
Let $E\rightarrow Y$ denote a holomorphic vector bundle on
a compact complex $D$-manifold $Y$. With notations as in 
Definitions \mb{\ref{Eulerdef}} and \mb{\ref{cherninv}},
$$
\chi(E) = \int_Y \mbox{Td}(Y) \mbox{ch}(E).
$$ 
\end{theorem}
To appreciate this theorem, recall that representatives of the Chern classes $c_k(Y)$, $c_k(E)$
can be obtained in terms of the curvature forms of arbitrary Hermitean metrics on $Y$ and $E$,
respectively (see e.g.\ \cite{lami89,bgv92} for a detailed discussion). 
The theorem thus yields the topological invariant $\chi(E)$ in terms of 
the integral of a local curvature expression. 
The essence of this theorem is the fact that 
$\slashed{D}^E=\partial_E+\overline\partial_E$
 is a Dirac operator  on $E=E^+\oplus E^-$, 
$\slashed{D}^E\colon\Gamma(E^\pm)\rightarrow\Gamma(E^\mp)$,
whose index 
$\dim\ker(\slashed{D}^E_{\mid E^+})-\dim\ker(\slashed{D}^E_{\mid E^-})=\chi(E)$ 
can therefore alternatively be calculated by heat kernel methods in terms of an integral
over local curvature data.

To recover the invariants of Definition \ref{Eulerdef},
for the holomorphic Euler characteristic $\chi(\mathcal O_Y)$
one uses the trivial bundle $E$,
so $\slashed{D}^E=\partial+\overline\partial$
and $\mbox{Td}(Y) \mbox{ch}(E)=\mbox{Td}(Y)$ with expansion
into Chern classes as in Definition \ref{cherninv}.
The usual Euler characteristic $\chi(Y)$ arises for
the virtual bundle $E=E^+-E^-$ with
$E^+=\oplus_{p\equiv0(2)} \Lambda^p T^\ast$, $E^-=\oplus_{p\equiv1(2)} \Lambda^p T^\ast$
and $T=T^{1,0}Y$ as before. With
$\chi(E) := \chi(E^+)-\chi(E^-)=\chi(Y)$, and 
 with notations as in Definition \ref{cherninv},
$$
\begin{array}{rcl}
\mbox{ch}(E) \;=\;\displaystyle \mbox{ch}(E^+)-\mbox{ch}(E^-)
&=&\displaystyle
 \prod_{j=1}^D (1-\exp(-y_j)),\;\\
 \mbox{Td}(Y) \mbox{ch}(E)
&=&\displaystyle \prod_{j=1}^D {y_j} \;=\;  c_D(Y).
\end{array}
$$
Thus the Hirzebruch-Riemann-Roch Formula \ref{ashrr}  yields
\begin{equation}\label{hrrappl}
\begin{array}{rcl}
\chi(\mathcal O_Y) 
&=&\displaystyle \int_Y \left( c_0(Y)+{\textstyle{1\over2}} c_1(Y) + {\textstyle{1\over12}} (c_1(Y)^2+c_2(Y)) 
+ \cdots \right), \\[1em]
\chi(Y) 
&=&\displaystyle \int_Y c(Y).
\end{array}
\end{equation}
To obtain the signature $\sigma(Y)$ from the Hirzebruch-Riemann-Roch Formula,
one uses the same total bundle $E^+\oplus E^-=\oplus_p \Lambda^p T^\ast$ with different $\Z_2$-grading
$\oplus_p \Lambda^p T^\ast=\tilde E^+\oplus\tilde E^-$ 
 and $E=\tilde E^+-\tilde E^-$.\\[0.5em]
There is another relation between the Hodge numbers
of $Y$ and its top Chern class, due to the interpretation
of $\chi(Y)$ in terms of the 
\textsc{Poincar\'e-Hopf Index Theorem}
\cite{po85a,ho26}. 
One combines this classical result from differential
topology, which is discussed in textbooks like \cite{gupo74,mo01},
with the celebrated Weitzenb\"ock formula from differential geometry, see 
e.g.\ the textbook \cite{jo95}.
Indeed, if $\chi(Y)\neq0$, then  the Poincar\'e-Hopf Index Theorem
implies that every holomorphic one-form on $Y$ has at least one zero.
On the other hand, using a Ricci-flat metric on $Y$, which
exists by the Calabi-Yau Theorem \ref{CYthm}, 
the Weitzenb\"ock formula implies that every holomorphic
one-form on $Y$ has constant norm. In other words:
\begin{lemma}\label{nonzeroEuler}
If a Calabi-Yau D-manifold $Y$ has non-vanishing Euler characteristic,
$\chi(Y)\neq0$, then $h^{1,0}(Y)=0$. 
\end{lemma}
We are now ready to classify Calabi-Yau two-manifolds; 
the resulting Theorem \ref{CYtwofolds} was first proved in \cite{ko64},
but with more recent results the proof can be simplified.
First note
\begin{lemma}\label{CY2Hodge}
Let $Y$ denote a connected Calabi-Yau two-manifold. Then $Y$ has Hodge numbers 
$$
h^{1,0}(Y) =2,\; h^{1,1}(Y) =4 \quad\mbox{ or }\quad
h^{1,0}(Y) =0,\; h^{1,1}(Y) =20.
$$
This fixes all Hodge numbers of $Y$.
\end{lemma}
\begin{proofsketch}
Since $D=2$, using  (\ref{hodgenumbers}) it is clear that
the values of
$h^{1,0}(Y)$ and $h^{1,1}(Y)$ determine all the Hodge numbers
of $Y$. 

Generalizing Lemma \ref{nonzeroEuler}, if a holomorphic $k$-form on $Y$
has a zero, then it vanishes identically on $Y$
 \cite{bo74}. This implies $h^{1,0}(Y)\leq2$. 
 Since $Y$ is K\"ahler, such that $h^{1,1}(Y)\neq0$, one
deduces that $\chi(Y)=0$ can only hold if 
$h^{1,0}(Y) =2,\; h^{1,1}(Y) =4$.

On the other hand, if $\chi(Y)\neq0$, then 
Lemma \ref{nonzeroEuler} implies $h^{1,0}(Y) =0$,
so by (\ref{hodgenumbers}) we have $h^{1,1}(Y)=\chi(Y)-4$.
Moreover, by (\ref{hodgenumbers}) and Definition \ref{Eulerdef}, 
$Y$ has holomorphic Euler characteristic $\chi(\mathcal O_Y)=2$. 
Hence   
$$
\chi(Y) 
\stackrel{(\ref{hrrappl})}{=} \int_Y c_2(Y) 
\stackrel{(\ref{hrrappl}),\, c_1(Y)=0}{=}  12\chi(\mathcal O_Y)= 24
$$
and therefore $h^{1,1}(Y)=\chi(Y)-4 =20$.
\end{proofsketch}

\noindent
If $\chi(Y)=0$ for a Calabi-Yau two-fold $Y$, then 
one can use the Calabi-Yau Theorem \ref{CYthm} and
a result by Bochner and Yano \cite{boya53} to prove
that $Y$ is a torus $\C^2/L$ with a lattice $L\subset\C^2$ of rank $4$. 
If $\chi(Y)\neq0$, then Lemma \ref{CY2Hodge} implies that $Y$ is a \textsc{K3 surface},
according to
\begin{definition}\label{K3def}
A connected, compact complex surface $Y$ with trivial ca\-nonical bundle 
and $b_{1}(Y) =\dim_\C H^{1}(Y,\C) =0$
is called a \textsc{K3 surface}.
\end{definition}
It was conjectured independently by Andreotti and Weil \cite{we58}
and proved by Siu \cite{si83} that all K3 surfaces are K\"ahler 
and thus Calabi-Yau\footnote{It is not clear to me who introduced
the name ``K3 surface'', and when. The standard explanation and first
mention in writing that I am aware of is due to Weil, who in
\cite{we58},  declares 
``Dans la seconde partie de mon rapport, il s'agit des vari\'et\'es K\"ahleriennes
dites K3, ainsi nomm\'ees en l'honneur de Kummer, Kodaira, K\"ahler et de
la belle montagne K2 au Cachemire.''
The explanation of course relies on the K\"ahler property for all K3 surfaces,
which at the time was only conjectural. It probably also relies on the fact that 
among mountaineers, K2 is often understood as the most challenging
summit to the very day, which however had been vanquished only
 a few years  before  Weil completed his report \cite{we58}, namely in 1954.}.
Using Lemma \ref{CY2Hodge}, I have thus summarized
a derivation of the first claim of 
\begin{theorem}[\cite{ko64}]\label{CYtwofolds}
If $Y$ is a connected Calabi-Yau two-manifold, then
$Y$ is either a complex two-torus or a K3 surface.
Viewed as real four-manifolds,
all complex two-tori are diffeomorphic to one another, 
and all K3 surfaces are diffeomorphic 
to one another.
\end{theorem}
In summary,
one main ingredient to the proof of the classification result 
Theorem \ref{CYtwofolds} is the Atiyah-Singer Index Theorem
in the form (\ref{hrrappl}), which for Calabi-Yau two-manifolds $Y$
implies $\chi(Y)=12\chi(\mathcal O_Y)$. 
For Calabi-Yau three-manifolds $Y$, the corresponding formula does 
not suffice to fix the topological type of $Y$. Indeed, the
problem of classifying all Calabi-Yau three-manifolds is wide open 
-- the naive observation that
there are more independent Chern classes to keep under control is
precisely the source of the problem. 
\subsection{The Kummer construction}\label{Kummer}
In the previous section, I stated that all K3 surfaces are diffeomorphic 
to one and the same real four-manifold $X$ \cite{ko64}. 
The choice of a complex structure and K\"ahler class on $X$, of course,
greatly influences its geometric properties, for example the symmetries
of the K3 surface.
The \textsc{Kummer construction}, which shall be discussed in the present
section, amounts to a special choice of complex structure and 
(degenerate) K\"ahler class, governed by the geometry of an underlying
complex two-torus $T_L$:
\begin{definition}\label{sgKummer}
Let $T_L=\C^2/L$ denote a complex two-torus, where $L\subset\C^2$ 
is a lattice of rank $4$, and $T_L$ carries the complex structure and
K\"ahler metric induced from the standard complex structure and 
Euclidean metric on  $\C^2$. This Calabi-Yau two-manifold enjoys
a biholomorphic isometry $\kappa\in\Aut(T_L)$ of order $2$ which is
induced by $z\mapsto-z$ on $\C^2$. 

The quotient $T/\Z_2$  of $T_L$ by the group $\{\id,\kappa\}\cong\Z_2$
is called the \textsc{singular Kummer surface} with underlying torus $T_L$.
\end{definition}
The singular Kummer surface $T_L/\Z_2$ is indeed singular:
We choose $\epsilon>0$ and denote by $B_\epsilon(0)\subset\C^3$ 
the open ball of radius $\epsilon$ with respect to the standard Euclidean
metric on $\C^3$. Then the map 
$$
\mathfrak u: \C^2\longrightarrow \C^3,\quad 
\mathfrak u(z^1,z^2):= \left((z^1)^2,\, (z^2)^2,\, z^1z^2\right)
$$
descends to an open neighborhood $U_\epsilon\subset T_L/\Z_2$ of $0$ in
the singular Kummer surface, where it is also denoted $\mathfrak u$, such that
$\mathfrak u$ bijectively maps $U_\epsilon$ to 
$$
\mathfrak u(U_\epsilon) = \{ u= (u^1,\,u^2,\,u^3)  \in B_\epsilon(0) \mid u^1u^2=(u^3)^2\}.
$$
The map $\mathfrak u$ is biholomorphic upon restriction to
$$
U_\epsilon\setminus\{0\} \stackrel{\mathfrak u}{\longrightarrow} \mathfrak u(U_\epsilon)\setminus\{0\}. 
$$
The equation $u^1u^2=(u^3)^2$ of $\mathfrak u(U_\epsilon)$ in $B_\epsilon(0)$
immediately shows that 
$\mathfrak u(U_\epsilon)$ is a double cone with an isolated singularity at $u=0$.
Let us define the \textsc{minimal resolution} of this singularity:
\begin{definition}\label{blowup}
With notations as above, the point $0\in U_\epsilon$, and equivalently
its image $0\in \mathfrak u(U_\epsilon)$, is called 
a \textsc{singularity of type $A_1$}. Let 
$\mathcal K_\epsilon:= \mathfrak u(\overline{U_\epsilon})\setminus\{0\}$, and
with homogeneous coordinates 
$v=(v^1\colon v^2\colon v^3)\in\P^2$, 
$$
W_\epsilon:=\left\{ (u,v)\in \mathcal K_\epsilon\times\P^2 \mid u^j v^k=u^kv^j 
\;\; \forall j,\,k\in\{1,2,3\}\right\}.
$$
Moreover, let $V_\epsilon$ denote the interior of the closure $\overline W_\epsilon$ 
of $W_\epsilon\subset\C^3\times\P^2$. Then  
$$
\sigma\colon V_\epsilon \longrightarrow \mathfrak u(U_\epsilon),\quad
\sigma(u,v):=u
$$ 
is called the \textsc{blow-up} of the singularity $0\in \mathfrak u(U_\epsilon)$ 
of type $A_1$, and 
$E:=\sigma^{-1}(0)$ is its \textsc{exceptional divisor}.
\end{definition}
If $(u,v)\in W_\epsilon$, then $u\neq0$, and the defining equations of $W_\epsilon$
imply $u^1u^2=(u^3)^2$ and $(v^1\colon v^2\colon v^3)=(u^1\colon u^2\colon u^3)$.
Denote this by  $v=[u]\in\P^2$ and observe that $v^1v^2=(v^3)^2$ follows.  
Then for $t\in\R$ close to $t=0$, say $t\in(0,\delta)$ for an appropriate
 $\delta>0$, the map
 $t\mapsto \gamma_u(t):= ( tu,[u])\in W_\epsilon$ yields a smooth curve with 
$\lim_{t\rightarrow0} \gamma_u(t) =(0,[u])\in E$. In fact,
$$
E=\{ (u,v)\in \C^3\times\P^2 \mid u=0,\, v^1v^2=(v^3)^2\}, 
$$
so the 
exceptional divisor $E$ is biholomorphic to $\P^1$,
$$
\P^1\stackrel{\cong}{\longrightarrow} E\quad
\mbox{ under }\quad (t^1\colon t^2)\mapsto \left(\;0,\; ((t^1)^2\colon (t^2)^2\colon t^1t^2)\;\right). 
$$
The most important properties of the blow-up of a singularity of type $A_1$
are summarized in the following
\begin{prop}\label{blowupproperties}
The resolution $\sigma\colon V_\epsilon\longrightarrow\mathfrak u(U_\epsilon)$ 
of the singularity $0\in U_\epsilon$ of type $A_1$
yields a smooth complex two-manifold $V_\epsilon$.
The restriction $\sigma_{\mid V_\epsilon\setminus E}\colon 
V_\epsilon\setminus E \longrightarrow \mathfrak u(U_\epsilon)\setminus\{0\}$
is biholomorphic, and the
exceptional divisor $E\subset V_\epsilon$ is biholomorphic
to $\P^1$. 
Moreover, $V_\epsilon$ has trivial canonical bundle.
\end{prop}
\begin{proofsketch}
The proof of Proposition \ref{blowupproperties}
can be performed by a direct calculation.

For example, smoothness of $V_\epsilon$
close to the point $p_1:=(0, (1\colon0\colon0))\in E$ can be checked 
in the chart $U_1:=\{ (u,v)\in V_\epsilon \mid v^1\neq0\}$. Indeed, 
$(u^1,v^3)\mapsto ((u^1, u^1(v^3)^2,u^1v^3),(1\colon (v^3)^2\colon v^3))$
yields a smooth parametrization of $U_1$ near $p_1$. Analogously, one obtains smoothness
everywhere, and holomorphicity of changes of coordinates is immediate,
as are the claims about biholomorphicity of $\sigma_{\mid V_\epsilon\setminus E}$ 
and $E\cong\P^1$ by what was said above.

The triviality of the canonical bundle of $V_\epsilon$
follows since $\eta_{(z^1,z^2)}=dz^1\wedge dz^2$ on $\C^2$ 
descends to a section of the canonical bundle of $V_\epsilon$ over $V_\epsilon\setminus E$. 
On $U_1$ and with respect to coordinates $(u^1,v^3)$ as above, we have
$$
\eta_{(u^1,v^3)} = {du^1\wedge du^3\over2u^1} 
= {\textstyle{1\over2}} du^1\wedge dv^3
$$
as long as $u^1\neq0$. However, this implies that $\eta$ has a holomorphic 
continuation to all of $U_1$, which never vanishes. By proceeding
similarly for other coordinate neighborhoods, $\eta$
can be continued to a nowhere vanishing global section of the canonical 
bundle of $V_\epsilon$, which thus is trivial.
\end{proofsketch}
The \textsc{Kummer construction} is now summarized in
\begin{theorem}[Kummer construction]\label{Kummerconstruction}
The singular Kummer surface $T_L/\Z_2$ of Definition \mb{\ref{sgKummer}}
has $16$ singularities of type $A_1$, situated at $L/2L\subset T_L/\Z_2$. 
The complex surface $X$ obtained by blowing up each
of these singularities is a K3 surface.
\end{theorem}
\begin{proofsketch}
One immediately checks that the singular points of $T_L/\Z_2$ are precisely the
images of the fixed points of $\Z_2$ in $T_L$ under the quotient $T_L\longrightarrow T_L/\Z_2$.
If $[z]$ denotes the image of $z\in\C^2$ under the natural projection $\C^2\longrightarrow T_L$, then
$[z]$ is fixed under $\Z_2$ if and only if
$[z]=[-z]$, that is, if and only if $2z\in L$. 
Since the lattice $L$ has rank $4$, we find that $L/2L\cong\mathbb F_2^4$ 
contains $16$ points.
That each of these singularities is of type $A_1$ is also immediate
-- one translates the coordinates induced from $\C^2$ on an
open neighborhood of $y\in L/2L\subset T_L$ by $-y$
and finds that a 
punctured neighborhood of the origin
is then biholomorphically mapped to $U_\epsilon\setminus\{0\}$ as in Definition \ref{blowup}.

Now $X$ is obtained from the singular Kummer surface by blowing up all
the singularities, that is, by replacing a neighborhood $U_\epsilon$ of each
of the singular points by a copy of the blow-up $V_\epsilon$. 
Since $U_\epsilon\setminus\{0\}$
is biholomorphic to $V_\epsilon\setminus\{0\}$
and according to Proposition \ref{blowupproperties},
$X$ is a smooth complex surface with trivial canonical bundle. Moreover,
$X$ is compact and connected by construction. I claim that $b_{1}(X)=0$. 
Indeed, 
$\kappa\in\Z_2$ acts by multiplication by $-1$ on $H^1(T_L,\C)$, 
such that none
of the classes in $H^1(T_L,\C)$ can descend  to $H^1(X,\C)$.
Furthermore, according to Proposition \ref{blowupproperties}, 
the exceptional divisor of each blow-up is biholomorphic to $\P^1$ 
with $h^{p,q}(\P^D)=\delta_{p,q}$,
and thus the blow-ups
only contribute to the cohomology of $X$ in even degree. 

In summary, $X$ is a connected, compact complex surface with trivial
canonical bundle and $b_{1}(X)=0$. 
By Definition \ref{K3def}, I have shown that $X$ is a K3 surface.
\end{proofsketch}

Note that the standard K\"ahler metric on $\C^2$, which is 
the Euclidean one, induces a \textsc{degenerate} 
K\"ahler metric on $X$. Indeed, the induced metric assigns vanishing
volume to each irreducible component $E\cong\P^1$ of the exceptional divisor
and thus does not correspond to a smooth, Riemannian metric on $X$. 
\begin{definition}\label{Kummerdef}
The K3 surface $X$ obtained from a complex torus $T_L=\C^2/L$ 
by the Kummer construction of Theorem \mb{\ref{Kummerconstruction}} is
called a \textsc{Kummer surface}. In other words, a Kummer surface
carries the complex structure and \(degenerate\) K\"ahler structure induced
from the respective structures inherited by $T_L$ from 
 the  standard ones on $\C^2$.
\end{definition}
By \cite{ni75}, every K3 surface $X$ which is obtained from some singular surface
by the minimal resolution
of $16$ distinct singularities of type $A_1$ 
is biholomorphic to a Kummer surface. In other words, there exists a 
rank $4$ lattice $L\subset\C^2$ such that $X$ is biholomorphic to the
minimal resolution of the singular Kummer surface $T_L/\Z_2$. There are many 
examples of singular quartic hypersurfaces in $\P^3$
with $16$ singularities of type $A_1$, e.g.\ the one defined by 
$$
\left\{ z= (z^0\colon\cdots\colon z^3)\in\P^3 \;\left|\;
\vphantom{\sum}
\smash{\sum_{j=0}^3} (z^j)^4 - 4 \smash{\prod_{j=0}^3} z^j\right.\right\}.
$$  
Such singular Kummer surfaces were first studied by Kummer \cite{kuges}. 
The construction has become a classical one by now and is described in detail, for
example, in \cite{bhpv84}.
\section{Conformal field theory}\label{CFT}
As argued in the introduction, the Calabi-Yau geometry
which the previous section was devoted to plays
a crucial role in string theory. Indeed, a
so-called \textsc{non-linear sigma model construction}  is
predicted to yield a \textsc{superconformal field theory}
as the world-sheet theory for superstrings in a Calabi-Yau target
geometry. 

Unfortunately, for a generic Calabi-Yau manifold it is still hopeless 
to attempt a non-linear sigma model construction 
explicitly. The only exceptions, in general, are the complex tori,
which carry flat metrics, such that non-linear sigma models
are only little more complicated than free field
theories. Furthermore, orbifold constructions, which can be viewed
as generalizations of
the Kummer construction of Section \ref{Kummer}, are reasonably
well understood. 

Non-linear sigma models on K3 surfaces, which will be discussed
in Section \ref{K3theories}, furnish
an intermediate case in that many \textsc{K3 theories} are accessible
through orbifold procedures, and in that the moduli space of these
theories has been determined globally under a few additional
assumptions \cite{se88,ce90,asmo94,nawe00}. 
That superconformal field theories 
allow an independent, mathematical approach is a crucial
ingredient to that result. 
The current section therefore gives an overview on aspects 
of conformal field theory related to the particular
Calabi-Yau geometries that the previous section was focused 
on, namely the non-linear sigma models on complex tori (Section \ref{toroidal})
and their $\Z_2$-orbifolds (Section \ref{Z2}). 
Due to restrictions of space and time, I content myself
with giving an overview and pointing to some relevant literature.
\subsection{Toroidal superconformal field theories}\label{toroidal}
In this treatise, I do not attempt to give a definition of 
conformal field theory (CFT), though hopefully these notes 
are useful also for the CFT-novice, in that they discuss some
of the basic ingredients to CFT. I have summarized my own view on  a
definition of CFT elsewhere, see e.g. \cite{we10}.
To make best use of the restricted amount of space, I refer
the reader to the recent review  \cite{we14} as
a companion paper to the present lecture notes. 
Indeed, while in \cite{we14}  the main emphasis lies on CFT aspects,
here I focus on the geometric point of view.
In particular,  for the notions of holomorphic fields, 
operator product expansions (OPEs) and normal ordered products,
following  \cite{lewi78,frka80,bo86,ka96,frbe04},
see \cite[Sect.~2.1]{we14}, and for a summary of a
definition of conformal and superconformal field theories, see
\cite[Sect.~2.2]{we14}. As in this reference, in what follows, I solely discuss
\textsc{two-dimensional Euclidean unitary conformal field theories}. All
superconformal field theories (SCFTs) are assumed to be \textsc{non-chiral} and to
enjoy \textsc{space-time supersymmetry} as well as 
\textsc{$N=(2,2)$ world-sheet supersymmetry}.
\\[0.5em]
An important ingredient to SCFT is the
representation theory of the \mbox{(super-)} Virasoro algebra. Indeed,
the \textsc{space of states} of any of our SCFTs
is a  $\Z_2\times\Z_2$-graded
complex vector space 
$$
\begin{array}{r@{\!}c@{\!}lr@{\!}c@{\!}l}
\H\;=\;\H^{NS}\oplus\H^R,\quad
\mbox{ where }\quad
\H^{NS}&\;\;=\;\;&\H^{NS}_b\oplus\H^{NS}_f,\quad&
\H^R&\;\;=\;\;& \H^R_b\oplus\H^R_f,\\[0.5em]
\H_b&:=&\H^{NS}_b\oplus\H^R_b,\quad&
\H_f&:=&\H^{NS}_f\oplus\H^R_f.
\end{array}
$$
The subspaces $\H^{NS}$ and $\H^R$ are called \textsc{Neveu-Schwarz (NS)}
and \textsc{Ramond (R) sectors}, while $\H_b$ and $\H_f$ are the subspaces
of \textsc{bosonic} and \textsc{fermionic} states, respectively.
The space $\H$ carries a representation 
of two copies of a super-Virasoro algebra at central charges
$(c,\overline c)$, where for all examples discussed here we
have $c=\overline c$. As usual, the
zero-modes of left- and right-moving Virasoro fields
and $U(1)$-currents are denoted $L_0,\,\overline L_0$ and
$J_0,\, \overline J_0$, respectively\footnote{As in \cite[Sect.~2.2]{we14}, 
I assume a \textsc{compactness condition} for our SCFTs,
namely that the operators $L_0,\,\overline L_0$,
$J_0,\, \overline J_0$ are simultaneously diagonalizable,
and that simultaneous eigenspaces of $L_0$ and $\overline L_0$
are finite dimensional.}. Then space-time supersymmetry
implies that $(-1)^F:=e^{\pi i(J_0-\overline J_0)}$ 
acts as identity operator $\id$
on $\H_b$ and as $-\id$ on $\H_f$. 

Every SCFT possesses a 
modular invariant \textsc{partition function}
$Z(\tau,z)$ in two complex variables $\tau,\, z\in\C$ with $\Im(\tau)>0$, 
$$
\mbox{with }
q:=e^{2\pi i\tau}, \, y:=e^{2\pi iz}\colon\quad
Z(\tau,z)
=\tr_{\H_b} 
\left( q^{L_0-c/24} y^{J_0}\overline q^{\overline L_0-\overline c/24} \overline y^{\overline J_0} \right). 
$$
Space-time supersymmetry is tantamount to a decomposition of the partition function
into four sectors,
\begin{eqnarray}
Z(\tau,z) 
&=& {\textstyle {1\over2}}\left( Z_{NS}(\tau,z)+Z_{\widetilde{NS}}(\tau,z)+
Z_{R}(\tau,z)+Z_{\widetilde{R}}(\tau,z)\right),\nonumber\\[1em]
\mbox{for } S\in\{NS,\,R\}\colon\!\!\!\!\!\!\!\!\!\!\!\!\!\!\!\!\nonumber\\
Z_{S}(\tau,z)
&:=& \tr_{\H^{S}} 
\left( q^{L_0-c/24} y^{J_0}\overline q^{\overline L_0-\overline c/24} \overline y^{\overline J_0} \right),\nonumber\\
&&\mbox{ with } 
Z_{R}(\tau,z)
\;=\; (q\overline q)^{c/24} (y\overline y)^{c/6} Z_{NS}(\tau,z+{\textstyle {\tau\over2}}),\label{spectralflows}\\[1em]
Z_{\widetilde{S}}(\tau,z)
&:=& \tr_{\H^{S}} 
\left( (-1)^{J_0-\overline J_0}
q^{L_0-c/24} y^{J_0}\overline q^{\overline L_0-\overline c/24} \overline y^{\overline J_0} \right)
\;=\;  Z_{S}(\tau,z+{\textstyle {1\over2}}).\nonumber
\end{eqnarray}
Equation (\ref{spectralflows}) reflects an isomorphism $\H^{NS}\cong\H^R$
of representations of the $N=(2,2)$ superconformal algebra, known as
\textsc{spectral flow}, which holds due to our assumption of space-time
supersymmetry.\\[0.5em]
The simplest example of such an SCFT
is a \textsc{toroidal $N=(2,2)$ superconformal field
theory}, see \cite[Def.~5]{we14} for a more detailed account. 
At central charges  $(c,\overline c)=(3D,3D)$, $D\in\N$, 
such a theory is
in particular characterized
by the fact that the chiral algebras on the left and on the right contain a
$\mathfrak u(1)^{2D}$-current algebra each, along with $D$ left-moving
and $D$ right-moving Dirac fermions, the superpartners of the 
$\mathfrak u(1)$-currents. The corresponding fields
$j^k(z),\; k\in\{1,\ldots,2D\}$ and 
$\psi_l^\pm(z),\; l\in\{1,\ldots,D\}$, whose only non-vanishing
OPEs are
\begin{equation}\label{torusfields}
\begin{array}{rcl}
\forall k,\,l\in\{1,\ldots,2D\}\colon\quad
j^k(z)j^l(w) &\sim&\displaystyle {\delta^{k,l}\over(z-w)^2},\\[1em]
\forall k,\,l\in\{1,\ldots,D\}\colon\quad
\psi^+_k(z)\psi^-_l(w) &\sim& \displaystyle{\delta^{k,l}\over(z-w)},
\end{array}
\end{equation}
along with their right-moving analogues 
are called the \textsc{basic fields}. They are
examples of so-called \textsc{free fields} (see e.g.\ \cite{ka96}).

The space of states $\H$ of a toroidal SCFT
arises as \textsc{Fock space representation} 
of the modes of these basic fields, that is, of the super-Lie algebra
generated by $\id$ and
$a_n^k,\, (\psi_l^\pm)_m\in\End(\H^S),\, 
k\in\{1,\ldots,2D\},\,l\in\{1,\ldots,D\},\,S\in\{NS, R\}$, where $n\in\Z$, and
on the NS-sector $\H^{NS}$, $m\in\Z+{1\over2}$, while on the R-sector $\H^R$, $m\in\Z$,
$$
\begin{array}{rclcl}
[a_m^j, a_n^k] &:=& a_m^j a_n^k -  a_n^k a_m^j &=& m\delta^{j,k}\delta_{m+n,0}\id,\\[0.5em]
\{\psi_m^j, \psi_n^k\} &:=& \psi_m^j \psi_n^k +  \psi_n^k \psi_m^j &=& \delta^{j,k}\delta_{m+n,0}\id
\quad\forall j,\,k,\; m,\,n,
\end{array}
$$
and analogously for the right-movers. All other (super-)commutators between the
$a_n^k,\, (\psi_l^\pm)_m$ and their right-moving analogues vanish. One then has
$\H^{NS}=\oplus_{\gamma\in\Gamma} \H_\gamma$ with 
$\Gamma$ the so-called \textsc{charge lattice},
and $\H_\gamma$ a lowest weight representation of the above
super-Lie algebra of modes, with lowest weight vector
$|\gamma\rangle \in \H_\gamma$ such that
\begin{eqnarray*}
\forall j\in\{1,\ldots,2D\}\colon\quad
a_0^j |\gamma\rangle &=& \gamma^j |\gamma\rangle,\\
\forall k\in\{1,\ldots,D\}\colon\;
a_m^j |\gamma\rangle &=& \psi_{m+{1\over2}}^k|\gamma\rangle = 0
\;\; \forall m\in\Z \mbox{ with } m<0,
\end{eqnarray*}
and analogously for the right-moving modes. Here, 
$\Gamma\subset\R^{2D,2D}$, by which I mean
$\R^{2D,2D}=\R^{2D}\oplus\R^{2D}$ and
$$
\forall \gamma\in\Gamma:\quad
\gamma=(\gamma_L,\gamma_R) \mbox{ with }
\gamma_L,\, \gamma_R\in\R^{2D},\;\;
\langle\gamma,\gamma\rangle = \gamma_L\cdot\gamma_L - \gamma_R\cdot \gamma_R
$$
with the standard Euclidean scalar product $\cdot$
on $\R^{2D}$ and $\gamma_L=(\gamma^1,\ldots,\gamma^{2D})^T$, 
$\gamma_R=(\gamma^{2D+1},\ldots,\gamma^{4D})^T$.
The representation $\H^R$ is obtained from $\H^{NS}$ by
spectral flow, as was mentioned in the discussion of equation
(\ref{spectralflows}). 
A counting argument then shows
$$
\tr_{\H_\gamma} 
\left( q^{L_0-c/24} y^{J_0}\overline q^{\overline L_0-\overline c/24} \overline y^{\overline J_0} \right)
= {q^{{1\over2}\gamma_L\cdot\gamma_L}\overline q^{{1\over2}\gamma_R\cdot\gamma_R}
\over \left| \eta(\tau) \right|^{4D}} 
\cdot { \left| {\vartheta_3(\tau,z)\over\eta(\tau) } \right|^{2D}}
$$
with the Dedekind eta function $\eta(\tau)$ and
the standard Jacobi theta  functions $\vartheta_k(\tau,z),\, k\in\{1,\ldots,4\}$. 
Here, the charge lattice $\Gamma$ is an even selfdual lattice of signature 
$(2D,2D)$, given in terms of 
an embedding into $\R^{2D,2D}$, which is specified by the
decomposition $\gamma=(\gamma_L,\gamma_R)$ for every $\gamma\in\Gamma$ 
as above. Generalizing
Kac's holomorphic lattice algebras \cite{ka96} to the non-holomorphic
case, in \cite{kaor00} the lattice vertex operator algebra 
corresponding to such an indefinite lattice is described. One has:
\begin{prop}[\cite{cent85,na86}]\label{formoftorustheory}
A toroidal superconformal field theory is uniquely characterized by its
charge lattice $\Gamma\subset\R^{2D,2D}$. 
For a theory with charge lattice $\Gamma$, setting 
$$
Z_\Gamma(\tau) 
:= \sum_{\gamma=(\gamma_L,\gamma_R)\in\Gamma} 
{q^{{1\over2}\gamma_L\cdot\gamma_L}\overline q^{{1\over2}\gamma_R\cdot\gamma_R}\over
|\eta(\tau)|^{4D}},
$$
the 
four sectors of the partition function $Z(\tau,z)$ are
$$
\begin{array}{rclrcl}
Z_{NS}(\tau,z) &=& Z_\Gamma(\tau) \cdot {\textstyle \left| {\vartheta_3(\tau,z)\over\eta(\tau) } \right|^{2D}},
\quad 
&Z_{\widetilde{NS}}(\tau,z) 
&=& Z_\Gamma(\tau) \cdot{\textstyle \left| {\vartheta_4(\tau,z)\over\eta(\tau) } \right|^{2D}},\\[0.5em]
Z_{R}(\tau,z) 
&=& Z_\Gamma(\tau)\cdot {\textstyle \left| {\vartheta_2(\tau,z)\over\eta(\tau) } \right|^{2D}},
&Z_{\widetilde{R}}(\tau,z) 
&=& Z_\Gamma(\tau)\cdot {\textstyle \left| {\vartheta_1(\tau,z)\over\eta(\tau) } \right|^{2D}},
\end{array}
$$
hence
$$
Z(\tau,z) = Z_\Gamma(\tau) \cdot Z_{f}(\tau,z)\quad \mbox{ with }\quad 
Z_{f}(\tau,z):=
\textstyle{1\over2} \sum\limits_{k=0}^4  \left| {\vartheta_k(\tau,z)\over\eta(\tau) } \right|^{2D}.
$$
\end{prop}
From Proposition \ref{formoftorustheory} one also reads
the moduli space of such theories \cite{cent85,na86},
see \cite[Thm.~1]{we14}.\\[0.5em]
As mentioned above, the complex tori are the unique Calabi-Yau manifolds
for which a non-linear sigma model construction can be performed directly and
explicitly; in fact, for a torus $T_L=\C^D/L$ with $L\subset\C^D$ a lattice of rank $2D$,
the non-linear sigma model construction yields a toroidal SCFT as discussed above.
Geometrically, each $\mathfrak u(1)$-current $j^k(z)$, $k\in\{1,\ldots,2D\}$,
as in (\ref{torusfields}),
arises as the field corresponding
to a parallel tangent vector field, given by a standard Euclidean coordinate
vector field in $\R^{2D}\cong\C^D$. The real and imaginary parts of 
its fermionic superpartners are
the fields corresponding to the dual cotangent vector fields. 
The resulting non-linear sigma model also depends on the choice of a \textsc{B-field},
which in this setting can  be described by a constant, skew-symmetric
endomorphism of $\R^{2D}$. Identifying $\R^{2D}\cong\C^D$ with its dual by means of
the standard Euclidean scalar product, the \textsc{dual lattice} of $L$
is given by
$$
L^\ast := \left\{ \alpha\in\R^{2D} \mid \alpha\cdot\lambda\in\Z \;\;\forall\lambda\in L \right\},
$$
and the charge lattice of the resulting toroidal SCFT
is 
$$
\Gamma 
= \left\{ {\textstyle{1\over\sqrt2}} (\mu-B\lambda+\lambda,\mu-B\lambda-\lambda) 
\mid\lambda\in L,\;\mu\in L^\ast \right\}.
$$
Vice versa, every toroidal SCFT allows
a \textsc{geometric interpretation} in terms of a non-linear sigma model,
see e.g.\ \cite{nawe00}:
\begin{definition}\label{geomint}
Consider a toroidal $N=(2,2)$ SCFT at 
central charges $(c,\overline c)=(3D,3D)$ with charge lattice
$\Gamma\subset\R^{2D,2D}$. 
A \textsc{geometric interpretation} of the theory is any
choice of complementary
$2D$-dimensional subspaces $Y,\,Y^0\subset\R^{2D,2D}$, such that
$\R^{2D,2D}=Y\oplus Y^0$, both $Y$ and $Y^0$ are null, and
both $\Gamma\cap Y$ and $\Gamma\cap Y^0$ are rank $2D$
lattices.
\end{definition}
The terminology deserves explanation:

\noindent
Given a geometric interpretation $\R^{2D,2D}=Y\oplus Y^0$ of a theory
with charge lattice $\Gamma\subset\R^{2D,2D}$, 
without loss of generality
we can set 
$$
\Gamma\cap Y=\left\{{\textstyle{1\over\sqrt2}}(\mu,\mu)\mid \mu\in L^\ast\right\},\;\;
\Gamma\cap Y^0 = \left\{ {\textstyle{1\over\sqrt2}} (\lambda-B\lambda,-\lambda-B\lambda) 
\mid\lambda\in L \right\}
$$ 
for some rank $2D$ lattice $L\subset\R^{2D}\cong\C^D$, $L^\ast\subset\R^{2D}$ its dual, and 
$B$ some skew-symmetric linear endomorphism of $\R^{2D}$. By 
Proposition \ref{formoftorustheory} and the explanations preceding 
Definition \ref{geomint}, our theory then agrees with a non-linear sigma model on
$T_L=\C^D/L$ with B-field $B$. 

Given such a geometric interpretation of a toroidal SCFT, 
it is now also clear 
how certain \textsc{geometric symmetries} of the torus $T_L$ may induce
symmetries of a toroidal SCFT on $T_L$:
If a  symmetry of $T_L$ is given in terms of 
$A\in\End(\C^D)$, then $A$ has to act as lattice automorphism of $L$.
Hence  $A\in O(2D)$
and $A$ also acts 
as lattice automorphism of $L^\ast$. If in addition, $AB=BA$,
then $A$ acts as lattice automorphism on $\Gamma$ which respects
the embedding $\Gamma\subset\R^{2D,2D}$. One checks that the
induced action on the space of states $\H$ yields a symmetry of the
toroidal SCFT.
\subsection{$\Z_2$-orbifold conformal field theories}\label{Z2}
In Section \ref{Kummer}, I presented the classical Kummer construction,
which yields a K3 surface from the much simpler complex two-torus 
$T_L=\C^2/L$ by \textsc{$\Z_2$-orbifolding}. 
The construction begins with the projection to $T_L/\Z_2$, where
$\Z_2$ is generated by the  symmetry $\kappa$ of $T_L$
which is induced from $-\id$ on $\C^2$.

By the discussion at the end of the previous section, $\kappa$ also acts as symmetry
on the non-linear sigma model constructed on $T_L$ with an arbitrary
B-field $B$. The present section is devoted to the
lift of this $\Z_2$-orbifold procedure to the level
of CFT, to construct a new superconformal field
theory from the much simpler toroidal one. This string theory 
procedure was inspired by the techniques that were 
originally developed in the context of \textsc{Monstrous Moonshine} 
\cite{flm84} for holomorphic vertex algebras,
and it can be carried out in much greater generality  \cite{dhvw85,dhvw86},
that is, in arbitrary dimensions and with more general orbifolding-groups.
Indeed, already the Kummer construction can be generalized
to obtain K3 surfaces by orbifolding an
underlying complex two-torus $T_L=\C^2/L$ with the 
appropriate symmetry
by groups $\Z_3,\, \Z_4$ or $\Z_6$, or even by certain
examples of non-Abelian groups, see \cite{we00} and
references therein. In the following I provide
an overview of the basic
ideas behind the orbifold procedure in CFT.
For the sake of brevity, 
I content myself with the discussion of $\Z_2$-orbifolds
as an instructive example.\\[0.5em]
If not stated otherwise, in the following I consider
a toroidal superconformal field theory 
at central charges $(c,\overline c)=(3D,3D)$ with geometric interpretation
on a complex torus $T_L=\C^D/L$ and with some B-field $B$
as discussed in Section \ref{toroidal}. As before, 
the space of states of the SCFT
is denoted by $\H=\H_b\oplus\H_f$, and $\kappa$ denotes
the symmetry of 
order $2$ of this theory which is induced by 
$\C^D\longrightarrow\C^D$, $z\mapsto-z$. 
It acts by multiplication by $-1$ on the basic fields (\ref{torusfields}) and on 
the charge lattice $\Gamma$, and therefore
as linear involution on the bosonic and fermionic spaces of states 
$\H_b$ and $\H_f$, respectively.

Let us investigate the $\kappa$-invariant subsector
$\H_b^{\Z_2}$ of the bosonic space of states. 
With notations as in Proposition \ref{formoftorustheory},
a counting argument shows
\begin{eqnarray*}
\tr_{\H_b^{\Z_2}}
\left( q^{L_0-c/24} y^{J_0}\overline q^{\overline L_0-\overline c/24} \overline y^{\overline J_0} \right)
&=& {\textstyle{1\over2}} \left( \vphantom{\textstyle{1\over2}}Z_\Gamma(\tau) \cdot Z_f(\tau,z)
+ Z_{-\id}(\tau,z)\right)\\
\mbox{with }\quad
Z_{-\id}(\tau,z) &=&
 \left|{2\eta(\tau)\over \vartheta_2(\tau,0)}\right|^{2D}\cdot Z_f(\tau,z).
\end{eqnarray*}
This expression is 
not modular invariant, since $Z_{-\id}(\tau,z)$ isn't, 
while $Z_\Gamma(\tau)$ and $Z_f(\tau,z)$ are. Hence 
$\H_b^{\Z_2}$ is not the space of bosonic states of a full-fledged superconformal
field theory. However, one can construct a \textsc{twisted sector} $\H_b^{tw}$ 
such that $\H_b^{\Z_2}\oplus\H_b^{tw}$ has this property.
To obtain information about such a would-be twisted sector,  we observe that 
$Z_{-\id}(\tau,z)$ has a natural modular invariant completion 
$$
\textstyle
Z_{-\id}\left(\tau,z\right)+Z_{-\id}\left(-{1\over\tau},{z\over\tau}\right)
+Z_{-\id}\left(-{1\over\tau+1},{z\over\tau+1}\right). 
$$
On the basis
of a path integral interpretation of each summand, see e.g.\
\cite[\S8.3]{gi88b},
modular invariance of this expression is expected; it is shown
to hold by a direct calculation.
The key idea behind the orbifolding procedure is to
interpret the three summands 
of the above expression as follows:
\begin{eqnarray*}
\textstyle
Z_{-\id}\left(\tau,z\right) 
&=& \textstyle
\tr_{\H_b}
\left( \kappa q^{L_0-c/24} y^{J_0}\overline q^{\overline L_0-\overline c/24} \overline y^{\overline J_0} \right),\\
\textstyle
Z_{-\id}\left(-{1\over\tau},{z\over\tau}\right)
&=& \textstyle
\tr_{\widetilde\H_b^{tw}}
\left( q^{L_0-c/24} y^{J_0}\overline q^{\overline L_0-\overline c/24} \overline y^{\overline J_0} \right),\\
\textstyle
Z_{-\id}\left(-{1\over\tau+1},{z\over\tau+1}\right) 
&=& \textstyle
\tr_{\widetilde\H_b^{tw}}
\left( \kappa q^{L_0-c/24} y^{J_0}\overline q^{\overline L_0-\overline c/24} \overline y^{\overline J_0} \right).
\end{eqnarray*}
Here, $\widetilde\H_b^{tw}$ needs to be constructed as
some representation  of the super-Virasoro algebras
that enjoys a symmetry $\kappa$ of order $2$ such that
$\H_b^{tw}$ is the $\kappa$-invariant subspace of $\widetilde\H_b^{tw}$. 
Indeed, in the present example 
an appropriate space $\widetilde\H_b^{tw}$ can be found, that 
is, such that $\H_b^{\Z_2}\oplus\H_b^{tw}$ can be interpreted as
bosonic space of states of a full-fledged superconformal field theory.
The details of the construction are quite technical but have been
worked out, see \cite{ffrs10} and references therein\footnote{Also
see Yi-Zhi Huang's blog \cite{hu14} and the references
therein to appreciate the mathematical problems that the
introduction of a SCFT on $\H_b^{\Z_2}\oplus\H_b^{tw}$ poses.}. 
\begin{prop}\label{Z2orbifoldsexist}
Consider a toroidal $N=(2,2)$ superconformal field theory 
at central charges $(c,\overline c)=(3D,3D)$ with 
charge lattice $\Gamma\subset\R^{2D,2D}$ and space of
states $\H$. 
Let $\kappa$ denote the symmetry of order $2$ which acts
by multiplication by $-1$ on $\Gamma$ and on the basic fields
\mb{(\ref{torusfields})}. 

Then there exists a \textsc{$\Z_2$-orbifold conformal field theory}
of this toroidal theory, whose space of states is $\H^{\Z_2}\oplus\H^{tw}$,
where $\H^{\Z_2}\subset\H$ denotes the $\kappa$-invariant subspace.
The four sectors of the partition function 
$Z^{orb}(\tau,z)$ of this theory are obtained
from its $\widetilde{R}$-sector by application of the spectral flow
formulas \mb{(\ref{spectralflows})}. With notations as in Proposition \mb{\ref{formoftorustheory}}
we have
\begin{eqnarray*}
Z^{orb}_{\widetilde R}(\tau,z)
&=& {\textstyle {1\over2}}\cdot\left(
Z_\Gamma(\tau)\cdot  \left| {\vartheta_1(\tau,z)\over\eta(\tau) } \right|^{2D}
\right.\\
&&\quad\quad\left.
+   \left| {2\vartheta_2(\tau,z)\over\vartheta_2(\tau,0) } \right|^{2D}
+ \left| {2\vartheta_4(\tau,z)\over\vartheta_4(\tau,0) } \right|^{2D}
+  \left| {2\vartheta_3(\tau,z)\over\vartheta_3(\tau,0) } \right|^{2D}
\right).
\end{eqnarray*}
\end{prop}
If the $\Z_2$-orbifold CFT of Proposition \ref{Z2orbifoldsexist}
is obtained from a toroidal SCFT with
geometric interpretation on $T_L=\C^D/L$ with some B-field
$B$, then it is believed that the orbifold theory can
be obtained by a non-linear sigma model construction from the
orbifold limit of the Calabi-Yau $D$-manifold obtained from
blowing up all singularities in $T_L/\Z_2$. 
To support this belief, one checks from the partition function 
$Z^{orb}(\tau, z)$ given in Proposition \ref{Z2orbifoldsexist} that the
dimension of the space of \textsc{twisted ground states}
in $\H^{tw}$ is $2^{2D}$. This is 
precisely the number of singular points in $T_L/\Z_2$. 
For more complicated orbifolding-groups, each singular 
point corresponds to a higher dimensional subspace
of twisted ground states, depending on the order
of the stabilizer group. This supports the idea that,
at least naively\footnote{Note however,
at least in the $D=2$-dimensional case, 
a subtlety concerning the B-field on the orbifold
\cite{nawe00}, which turns out to be non-zero 
on every irreducible component of the exceptional divisor
in the blow-up.}, the
introduction of twisted sectors corresponds
to the resolution of the singular points in the non-linear
sigma model on $T_L/\Z_2$. 

If $D=2$, then by Theorem \ref{Kummerconstruction}
the resolution of all singularities in $T_L/\Z_2$ yields
a K3 surface. We should
hence expect the $\Z_2$-orbifold conformal field theory 
of Proposition \ref{Z2orbifoldsexist} to allow a 
non-linear sigma model interpretation on a Kummer 
surface in this case. 
In \cite{nawe00}, a map between
the respective moduli spaces of conformal field theories
is constructed which is compatible with this expectation. 
Earlier evidence in favor of the prediction 
arises from a calculation
of the so-called \textsc{conformal field theoretic elliptic genus} 
\cite{eoty89}, according to
\begin{definition}\label{CFTellgen}
Consider the $\widetilde R$-sector $Z_{\widetilde{R}}$
of the partition function of an $N=(2,2)$-superconformal field theory,
viewed as a function of four complex variables 
$Z_{\widetilde{R}}=Z_{\widetilde{R}}(\tau,z;\overline\tau,\overline z)$. 
The \textsc{conformal field theoretic elliptic genus} of the theory is
$$
\mathcal E(\tau,z) 
:= Z_{\widetilde{R}}(\tau,z;\overline\tau,\overline z=0).
$$
\end{definition}
For a CFT that arises as 
non-linear sigma model on some Calabi-Yau $D$-manifold $Y$, it is expected
that the conformal field theoretic elliptic genus agrees with
the \textsc{geometric elliptic genus} of $Y$.
The latter can be defined as the
holomorphic Euler characteristic of
a formal vector bundle $\mathbb E_{q,-y}$ on $Y$,
more precisely a formal power series in $q=e^{2\pi i\tau}$ 
and $y=e^{2\pi iz}$ whose coefficients are holomorphic vector bundles on $Y$.
Using the notations $\Lambda_x E$, $S_x E$
for any vector bundle 
$E\rightarrow Y$ and a formal variable $x$ with 
$$
\Lambda_x E := \bigoplus_{p=0}^\infty x^p \Lambda^p E,\quad
S_x E:= \bigoplus_{p=0}^\infty x^p S^p E,
$$
where $\Lambda^p E,\, S^p E$ denote the $p^{\rm th}$ exterior 
and symmetric powers of $E$, and with  $T:=T^{1,0}Y$ the
holomorphic tangent bundle of $Y$,
\begin{equation}\label{bdlforellgen}
\mathbb E_{q,-y}  
= y^{-D/2}  \bigotimes_{n=1}^\infty \left( \Lambda_{-yq^{n-1}} T^\ast\otimes \Lambda_{-y^{-1}q^n} T
\otimes S_{q^n} T^\ast\otimes S_{q^n} T\right).
\end{equation}
The holomorphic Euler characteristic of Definition \ref{Eulerdef} is  naturally extended to  
formal power series with coefficients in holomorphic
vector bundles on $Y$,
$$
\chi\left( \smash{\sum_{n=0}^\infty} \vphantom{\sum}
x^n F_n\right):=\sum_{n=0}^\infty x^n \chi\left(F_n\right). 
$$
We then have:
\begin{definition}\label{geomellgen}
For a Calabi-Yau $D$-manifold $Y$ with holomorphic tangent
bundle $T:=T^{1,0}Y$, the \textsc{geometric elliptic genus}
$\mathcal E_Y(\tau,z)$ is the holomorphic Euler characteristic
of the bundle $\mathbb E_{q,-y}$ introduced in \mb{(\ref{bdlforellgen})},
$$
\mathcal E_Y(\tau,z):= \chi( \mathbb E_{q,-y} ).
$$
\end{definition}
Following \cite{wi87}, the geometric elliptic genus can
be interpreted as a regularized version of a $U(1)$-equivariant
index of a Dirac operator on the loop space of $Y$.
It is a topological invariant of $Y$ with many beautiful
mathematical properties. In particular, its modular transformation
properties agree with those of the conformal field theoretic elliptic
genus of Definition \ref{CFTellgen} for a theory at central charges
$c=\overline c=3D,\, D\in\N$: Both are \textsc{weak Jacobi forms} of weight $0$ and
index ${D\over2}$ (with a character, if $D$ is odd) 
\cite{hi88,wi88,eoty89,kr90,dfya93,wi94,boli00}.
See Section \ref{K3theories} and
\cite[Sect.~2.4]{we14} for a more detailed discussion of the 
two versions (geometric vs.\ conformal field theoretic)
of the elliptic genus and their expected
relationship, and for further references.\\[0.5em]
From Propositions \ref{formoftorustheory} and  \ref{Z2orbifoldsexist} one 
immediately finds the conformal field theoretic elliptic genera
of the non-linear sigma model
on a complex two-torus $T_L=\C^2/L$ and
of its $\Z_2$-orbifold CFT:
\begin{equation}\label{K3ellgen}
\begin{array}{rcl}
\mathcal E_{T_L}(\tau,z)
&=&0,\\
\mathcal E_{K3}(\tau,z)
&=&\displaystyle  8 \left( {\vartheta_2(\tau,z)\over\vartheta_2(\tau,0)}\right)^{2}
+8 \left( {\vartheta_3(\tau,z)\over\vartheta_3(\tau,0)}\right)^{2}
+8 \left( {\vartheta_4(\tau,z)\over\vartheta_4(\tau,0)}\right)^{2}.
\end{array}
\end{equation}
These functions agree with the known geometric elliptic genera
of $T_L$ and a Kummer surface, respectively, and thereby of all 
complex two-tori resp.\ K3 surfaces \cite{eoty89}:
\begin{prop}\label{Z2orbifoldisK3theory}
Consider a toroidal $N=(2,2)$ superconformal field theory at 
central charges $c=\overline c=6$ with geometric interpretation
on a complex two-torus $T_L=\C^2/L$ with some B-field $B$.
Then its conformal field theoretic elliptic genus 
agrees with the geometric elliptic genus of $T_L$, 
and thereby of every complex two-torus.
The conformal field theoretic elliptic genus of 
the $\Z_2$-orbifold CFT of this toroidal theory 
agrees with the geometric elliptic genus of any 
\(and thereby every\) K3 surface.
\end{prop}
\section{Outlook: Towards superconformal field the\-o\-ry on K3 and beyond}\label{K3theories}
\setcounter{definition}{0}
\renewcommand{\thedefinition}{\thesection.\arabic{definition}}
The discussion in the previous section foreshadows a rather
subtle relation between geometry and conformal field theory.
This final section of these lecture notes gives
a rough overview and outlook on 
attempts to get a better understanding of this relation.
Special attention is paid to K3 surfaces and conformal
field theories associated to them,
because as we have seen, the K3 surfaces furnish the
simplest examples where non-linear sigma model
constructions are not fully understood.\\[0.5em]
From a mathematical point of view, an abstract approach to 
CFT, based in representation theory, is
desirable. However, then the
route back to geometry is not immediate\footnote{Specifically concerning the
discussion of the ``geometry of CFTs'', see also
 Yi-Zhi Huang's blog \cite{hu14} and the references
therein.}. As mentioned before,
direct constructions of superconformal field theories
from Calabi-Yau data are sparse -- they are 
essentially restricted to toroidal SCFTs
and their orbifolds.

However, if the existence of a non-linear sigma model 
on a Calabi-Yau $D$-manifold $Y$ is assumed, then a number of 
additional properties  are known for the resulting SCFT.
First, it enjoys $N=(2,2)$ (world-sheet) supersymmetry at
central charges $c=\overline c=3D$ as well as space-time supersymmetry.
Second, all eigenvalues of the linear operators
$J_0$ and $\overline J_0$ on the space of states $\H=\H^{NS}\oplus\H^R$ are integral in 
the Neveu-Schwarz sector\footnote{Since space-time supersymmetry
implies  an equivalence of representations of the $N=(2,2)$ superconformal algebra
$\H^{NS}\cong\H^R$ 
under spectral flow, 
this integrality condition implies that 
$J_0-\overline J_0$ has only integral eigenvalues on all of $\H$,
and that the eigenvalues of $J_0$ and $\overline J_0$ in the
Ramond sector $\H^R$  lie in ${D\over2}+\Z$.} $\H^{NS}$.

At small central charges, more precisely at $c=\overline c=3$ and
$c=\overline c=6$, these conditions severely restrict the 
types of theories that can occur.
One first observes that
the conformal field theoretic elliptic genus 
of Definition \ref{CFTellgen} either vanishes,
or it agrees with the geometric elliptic genus
$\mathcal E_{K3}(\tau,z)$  of K3 surfaces\footnote{This can
only happen if $c=\overline c=6$.} as in (\ref{K3ellgen}):
For $c=\overline c=3$ this follows immediately, since the space
of weak Jacobi forms of weight $0$ and index ${1\over2}$ is trivial,
for $c=\overline c=6$ see \cite[Prop.~2(1)]{we14}.
This reproduces
the classification of Calabi-Yau $D$-manifolds with $D\leq2$
according to
Theorems \ref{CYonefolds} and \ref{CYtwofolds}
on the level of the elliptic genus. 
If the conformal field theoretic elliptic genus vanishes
for a theory at $c=\overline c=3$ or $c=\overline c=6$
which obeys all the above assumptions on supersymmetry
and on the $J_0,\,\overline J_0$-eigenvalues, 
then one actually finds that the
underlying theory is toroidal (see \cite[Prop.~2(2)]{we14}
for the case $c=\overline c=6$ -- the case $c=\overline c=3$ 
is in fact simpler and can be treated analogously).
This motivates the following definition, see also \cite[Def.~8]{we14}:
\begin{definition}\label{K3theoriesdef}
A superconformal field theory is called a \textsc{K3 theory}, if the following
conditions hold: The CFT is an $N=(2,2)$ superconformal field theory at
 central charges $c=\overline c=6$ with space-time supersymmetry, 
all the eigenvalues of $J_0$ and of $\overline J_0$ are integral,
and the conformal field theoretic elliptic genus of the theory agrees
with the geometric elliptic genus of K3 surfaces 
$\mathcal E_{\rm K3}(\tau, z)$
as in \mb{(\ref{K3ellgen})}.
\end{definition}
By Proposition \ref{Z2orbifoldisK3theory}, all
$\Z_2$-orbifold conformal field theories obtained from toroidal
superconformal ones at  $c=\overline c=6$ are examples of K3 theories.\\[0.5em]
One needs to appreciate that
Definition \ref{K3theoriesdef} does not make use of non-linear 
sigma model assumptions. This also
means that one needs
to carefully distinguish between \textsc{K3 theories} and 
\textsc{conformal field theories on K3}. It is wide open,
yet interesting and important, whether all K3 theories 
are theories on K3 in the sense that they can be constructed
as non-linear sigma models. At least to our
knowledge, no counter example is known.

There are various geometric properties of K3 surfaces
that can be recovered
from abstractly defined K3 theories. Let us only mention
that under few additional assumptions, listed in
\cite{se88,ce90,asmo94,nawe00},
each connected component of
the moduli space of K3 theories agrees with the
moduli space expected for non-linear sigma models on K3. 
From this identification one obtains
 the notion of \textsc{geometric
interpretation} for K3 theories in the spirit of 
Definition \ref{geomint}, see \cite{asmo94}.
For  the so-called \textsc{Gepner models},
this allows to 
interpret certain symmetries of these theories in terms of
geometric symmetries of an underlying K3 surface, see e.g.\ 
\cite{eoty89,nawe00,we01,we03,we05}.
The precise role of geometric symmetries for the properties
of the elliptic genus is  central to more recent discussions
\cite{tawe11,tawe12,tawe13}
of the \textsc{Mathieu Moonshine} phenomena, in search of 
a geometric explanation for the seminal observations of \cite{eot10,ga12},
see \cite[Sect.~4]{we14} for a summary.
Using mirror symmetry, one can also explicitly
construct a family of K3 theories with a geometric
interpretation on a family of smooth K3 surfaces, which
are expected to provide examples of 
non-linear sigma models on such smooth manifolds
\cite{we05}.\\[0.5em]
It is natural, but much more subtle, to 
reconstruct aspects of the vertex algebra which underlies
a superconformal field theory from geometric data,
if the theory is expected to have a non-linear
sigma model interpretation.
From a string theory perspective, one should 
introduce free fields
on local coordinate patches of $Y$, analogously
to the fields $j^k(z),\, \psi_l^\pm(z)$ of (\ref{torusfields}),
however it is unclear how to relate these ``locally defined'' 
fields to the ones that are used in abstract conformal field
theories. 

A more ``rigid'' approach, inspired by and closely related
to \textsc{topological field theory},
uses $bc-\beta\gamma$ systems on holomorphic
coordinate patches of a compact, connected complex
manifold $Y$. Roughly, one replaces
the local holomorphic
coordinate functions and holomorphic tangent vector fields
by free bosonic fields of scaling dimensions zero  and one,
respectively, (the \textsc{$bc$-fields}),
while cotangent vector fields correspond to 
fermionic fields (the \textsc{$\beta\gamma$-fields}).
For these fields, 
transition functions between holomorphic coordinate patches 
can be defined according to
the known geometric transformation rules. By this 
construction, 
one arrives at a sheaf of conformal vertex algebras on $Y$,
which allows the definition of a natural $\Z$-grading and 
a differential of degree one. The resulting complex is
known as the 
\textsc{chiral de Rham complex} \cite{msv98} because
it  is  quasi-isomorphic to the classical
(Dolbeault-) de Rham complex on  $Y$. 

The construction of the chiral de Rham complex
is ``rigid'' or topological in that an underlying conformal
 covariance manifests itself in a global 
section of the chiral de Rham complex, which yields
a Virasoro field at central charge zero. 
One can perform a topological twist to a traditional $N=2$ superconformal
structure at central charge $c=3D$ if $Y$ is a Calabi-Yau 
$D$-manifold. 
This structure is  found to
descend to the sheaf cohomology
of the chiral de Rham complex, which  carries
the structure of a superconformal vertex operator algebra 
\cite{bo01,boli00}. Thereby, the above-mentioned problem
is solved,
that the geometrically motivated fields are only defined locally
on coordinate patches, while abstractly defined conformal
field theories do not exhibit such ``locally defined'' fields.

Armed with this success, one would
naturally expect the vertex algebra that is obtained from
the chiral de Rham complex
to be closely related to
the one which underlies a non-linear sigma model on $Y$. 
This expectation is reinforced by the fact that almost by construction,
the \textsc{graded Euler characteristic} of the chiral de Rham
complex on a Calabi-Yau $D$-manifold $Y$
yields the geometric elliptic genus $\mathcal E_Y(\tau,z)$ 
of Definition \ref{geomellgen} \cite{bo01,boli00}. 
However, the
precise relation between the relevant vertex algebras turns
out to be more subtle. In particular, the construction of the 
chiral de Rham complex does not depend on the  choice
of the K\"ahler class on $Y$, in contrast to what one expects
for the vertex algebras obtained in non-linear sigma models. 
According to \cite{ka05}, this problem can be solved by performing
a large volume limit to identify the BRST-cohomology of
a topologically half twisted non-linear sigma model on $Y$ with 
the sheaf cohomology of the chiral de Rham complex.\\[0.5em]
Returning to K3 theories, where the notion of geometric 
interpretations is well understood, one may hope for more
concrete results. Indeed, for generic Calabi-Yau $D$-manifolds
it is notoriously hard to calculate the sheaf cohomology of the
chiral de Rham complex and its superconformal vertex algebra
structure. On the other hand, 
for toroidal SCFTs and their $\Z_2$-orbifolds
at $c=\overline c=6$, I expect that
a direct calculation should show that the 
vertex algebra obtained from the  
sheaf cohomology of the chiral de Rham complex in fact agrees with
the one obtained by a topological half-twist from the respective 
superconformal field theories \cite{cgw15}. Intriguingly, our
calculations seem to indicate that 
although the construction of the chiral de Rham complex crucially
depends on the choice of complex structure on $Y$, the
resulting vertex algebras 
show no dependence on the choice of complex structure
(nor K\"ahler structure)
within these two classes of examples. \\[0.5em]
In short, there are many intriguing mathematical and physical
mysteries left
when it comes to  
the journey from geometry  to
conformal field theory, even en route at the potentially 
simpler K3 geometry and K3 theories. 
%
%
%
%
%
\bibliographystyle{kw}
\bibliography{kw} 
 \end{document}